\documentclass[10pt]{article}

\usepackage{amssymb}

\newtheorem{theorem}{Theorem}
\newtheorem{definition}[theorem]{Definition}
\newtheorem{lemma}[theorem]{Lemma}
\newtheorem{proposition}[theorem]{Proposition}
\newtheorem{corollary}[theorem]{Corollary}

\usepackage{showlabels}

\newcommand{\var}{\mbox{\rm var}}
\newcommand{\const}{\mbox{\rm const.}}
\def\qed{\unskip\nobreak\hfill\penalty50\hskip 3pt\hbox{}\nobreak
\hfill\hbox{\vrule width 4 pt height 10 pt}}

\begin{document}
\bibliographystyle{plain}

\title {The compound Poisson distribution and return times in dynamical systems}
\author{Nicolai Haydn\thanks{Mathematics Department, University of Southern California,
Los Angeles, 90089-1113. E-mail: $<$nhaydn@math.usc.edu$>$. This work was supported by a grant
from the NSF (DMS-0301910) and by a grant from the Universit\'e de Toulon et du Var, France.} 
\and Sandro
Vaienti\thanks{Centre de Physique Th\'eorique, UMR~6207,  CNRS, Luminy Case~907, F-13288
Marseille Cedex~9, and Universities
of Aix-Marseille~I,~II and Toulon-Var. F\'ed\'eration de Recherche des Unit\'es de
Math\'ematiques de Marseille. E-mail: $<$vaienti@cpt.univ-mrs.fr$>$.}}

\maketitle

\begin{abstract}
Previously it has been shown that some classes of mixing dynamical systems have limiting return
times distributions that are almost everywhere Poissonian. Here we study the behaviour of return
times at periodic points and show that the limiting distribution is a compound Poissonian
distribution. We also derive error terms for the convergence to the limiting distribution.
We also prove a very general theorem that can be used to establish compound 
Poisson distributions in many other settings.
\end{abstract}

\section{Introduction}

In 1899 Poincar\'e showed that for a map $T$ on some space $\Omega$ which has an 
invariant probability measure $\mu$, almost every every point returns within finite time 
arbitrarily close. This means that for every (measurable) $A\subset\Omega$ with $\mu(A)>0$
the return time function $\tau_A(x)=\min\{k\ge1:\;\; T^kx\in A\}$ is finite for $\mu$-almost every
$x\in A$. This result was quantified by Kac in 1947 for ergodic measures. His theorem 
states that $\int_A\tau_A(x)\,d\mu(x)=1$, provided $\mu$ is ergodic, which implies that
$\tau_A(x)$ is on average equal to $1/\mu(A)$. 
Since 1990 there has been a growing interest in the statistics of return times and 
in particular in the distribution of $\tau_A$. Considering that it was shown in \cite{Lac,KL}
that for ergodic measures the limiting distribution of a sequence of (rescaled)
return functions $\tau_{U_n}$ for $n\rightarrow\infty$
can be any arbitrarily prescribed distribution for suitably chosen sets $U_n$, 
it is necessary to assume that the return sets $A$ are dynamically regular. 

For the measure of maximal entropy on a subshift of finite type,  Pitskel \cite{Pitskel} 
showed that the return times are in the limit Poisson distributed for cylinder sets
$A_n(x)$ ($\mu(A_n(x))\rightarrow 0$ as $n\rightarrow\infty$) 
where the set of suitable `centres' $x\in\Omega$ form a full measure set. 
A similar result had independently been obtained by Hirata \cite{Hirata1,Hirata2}
by a different method. For $\phi$-mixing Gibbs measures Galves and Schmitt \cite{GS}
showed in 1997 that  the first return time is in the limit exponentially distributed and 
that the convergence is at an exponential rate. Subsequent results (e.g.\
\cite{WTW,Den,HSV,HV,Coe,AG}) established limiting distribution results for 
first or multiple return times in various settings and sometimes with rates of 
convergence. 

Almost all previous results look at the distribution of return times near generic
points. The notable exception being the paper \cite{Hirata1} by Hirata which gives the 
distribution of the first return time at a periodic point. In the present paper we consider 
periodic points for sufficiently well mixing invariant measures and show that
the limiting distribution is compound Poissonian. The compound Poisson distribution
 has previously
been used in various settings including the analysis of internet traffic where the 
waiting time between packets is exponential and the size of each packet is
geometrically distributed. It has also been used to model the survival of 
capitalist enterprises in the free market system~\cite{Med}.
 The main technical result, Proposition~\ref{sevastyanov}, provides
conditions under which one obtains a compound Poissonian distribution with 
error terms. This result in itself will be of interest to a much larger community than
the one of dynamicists addressed in this paper. 

\vspace{3mm}
\noindent
Let $\mu$ be a probability measure on a space $\Omega$ which carries a transformation
$T$, preserving $\mu$, and whose $\sigma$-algebra is generated by the joins ${\cal
A}^k=\bigvee_{j=0}^{k-1}T^{-j}{\cal A}$, $k=1,2,\dots$, of a given finite measurable partition $\cal A$.
The elements of ${\cal A}^k$ are called {\em $k$-cylinders}.
We assume $\cal A$ is generating, i.e.\ the elements of ${\cal A}^\infty$ are single points.
Denote
by $\chi_A$ the characteristic function of a (measurable) set $A$ and define the random variable:
$$
\zeta_A(z)=\sum_{j=1}^{\tau}\chi_A\circ T^j(z),
$$
$z\in \Omega$.
The value of $\zeta_A$ measures the number of times a given point returns to $A$ within the time
$\tau$. Typically the obsevation time $\tau $ is chosen to be $[t/\mu(A)]$ where $t$ is a
parameter value. (The rescaling factor $1/\mu(A)$ agrees with Kac's theorem.) For instance,
if $\mu$ is the measure of maximal entropy on a subshift of finite type, then
Pitskel~\cite{Pitskel} showed that $\zeta_{A_n(x)}$ is for $\mu$-almost every $x\in\Omega$
in the limit $n\rightarrow\infty$ Poisson distributed (where $\tau_n=[t/\mu(A_n(x))]$
and $A_n(x)$ denotes the unique $n$-cylinder that contains $x$).
In \cite{HV} we have proven a similar result for a much wider class  of systems and provided
error estimates.

We develop a mechanism which allows to prove the compound Poisson distribution of return times at
periodic points $x$ and also to obtain error estimates as the cylinder sets $A_n(x)$ shrink in
measure to zero.

To be more precise, if we denote by $\zeta_n^t(z)$  the counting function
$\sum_{j=1}^{\tau_n} \chi_{A_n(x)}\circ T^j(z)$, with the {\em observation time}
$\tau_n=\left[\frac{t}{(1-p)\mu(A_n(x))}\right]$, we will study the following distribution:
\begin{equation}\label{DE}
\mathbb{P}\left(\zeta_n^t=r\right),\;\;\;\; r=0,1,2,\dots,
\end{equation}
where $t>0$ is a parameter and $p\in[0,1)$
depends on the periodic point $x$ and will be given in Sect.~3.
We will show that the limit $n\rightarrow \infty$ is the compound Poisson
distribution (see Section~2) if $\mu$ is a $(\phi,f)$-mixing measure. We also provide
rates of convergence. This then implies under some mild additional assumptions~\cite{HLV} 
the uniform integrability of the process $\zeta_n^t$.

We then extend this result to return times, i.e.\ to the distribution of 
 the process $\zeta_n^t(z)$ restricted to the cylinder $A_n(x)$. 
The measure $\mu$ is then
replaced by the conditional measure $\mu_n=\left.\frac1{\mu(A_n(x))}\mu\right|_{A_n(x)}$. 
We refer to this
second case as the distribution of the number of visits for return times.

Our results for return times
considerably improve on the work of Hirata~\cite{Hirata1}, where he computed (without error)
the distribution of the first return time (order $r=0$) around periodic points and for Gibbs measures on
Axiom-A systems.

\vspace{3mm}

\noindent
The plan of the paper is the following.
 The purpose of  section~2 is to prove Proposition~\ref{sevastyanov} that gives general
conditions under which a sum of mutually dependent $0,1$-valued random variables converges to the
compound Poisson distribution and provides error terms. A similar result that had
been inspired by a theorem of Sevast'yanov~\cite{Sev}, was
proved and used  in \cite{HV} for the Poisson distribution.

The distribution of return times is tied to the mixing properties of the invariant measure
considered. For that purpose we introduce in the third section the $(\phi,f)$-mixing property. This
property is more general that the widely used $\phi$-mixing property and is reminiscent of Philipp
and Stout~\cite{PS} `retarded strong mixing property'. In this way one can obtain distribution
results on return times of some well studied dynamical systems that are not $\phi$-mixing, e.g.\
rational maps, parabolic maps, piecewise expanding maps in higher dimension~\dots.

The third
section is devoted to the proof of the existence of the limit distribution and rates of 
convergence for
entry times (Theorem~7), while the fourth section extends those results to return times 
(Theorem~10). 
Section~5 contains a careful application to rational maps with critical points (Theorem~11).

We conclude this introduction with an interesting observation. Limit distributions for entry and
return times have been provided along nested sequences of cylinder sets converging to points $x$
which were chosen almost everywhere or as periodic points. In section~3.4 we will show how to find
points $x$ which do not have limit distributions at all, and this will be achieved by using our
results on the compound Poisson distribution around periodic points.

\section{Factorial moments and mixing}
The main purpose of this section is to prove a very general result
which we use to prove the main results in sections~3 and~5 but which can also be useful to establish
compound Poisson distribution with respect to the geometric distributiion in many other settings. 
For more general compound Poisson distributions see~\cite{Fel}. More recently
(e.g.~\cite{CR,BCL}) there have
been efforts to approach compound Poisson distributions using the Chen-Stein method.
The treatment in~\cite{CR} has a more general setting, but the result is far from applicable to 
our situation.\footnote{We thank the referee for pointing us towards Chen and Roos' work and also for
other enlightening remarks.}
Proposition~\ref{sevastyanov} is of general interest and is reminiscent of 
existing theorems which establish the Poisson distribution (cf.\ \cite{Sev,HV}). 
from the convergence of the moments. 
It provides general conditions under which the distribution of a finite set of $0,1$-valued
random variables is close to compound Poisson (and provides error terms).
In  sections~3 and~5 we then use it to obtain the speed of convergence for the limiting
distributions for $\phi$-mixing systems, some non-Markovian systems and equilibrium states for
rational maps with critical points.

\subsection{Compound Poisson distribution}

For a parameter $p\in[0,1)$ let us define the polynomials
$$
P_r(t,p)=\sum_{j=1}^rp^{r-j}(1-p)^j\frac{t^j}{j!} \left(\begin{array}{c}r-1\\j-1\end{array}\right),
$$
$r=1,2,\dots$, where $P_0=1$ ($r=0$). The distribution $e^{-t}P_r(t,p)$, $r=0,1,2,\dots$ 
is  sometimes called the P\'olya-Aeppli distribution~\cite{JKW}.
It has the generating function
$$
g_p(z)=e^{-t}\sum_{r=0}^\infty z^rP_r=e^{t\frac{z-1}{1-pz}},
$$
its mean is $\frac{t}{1-p}$ and its variance is $t\frac{1+p}{(1-p)^2}$. 
Note
that for $p=0$ we recover the Poisson terms $e^{-t}P_r(t,0)=e^{-t}\frac{t^r}{r!}$ and the generating
function $g_0(z)=e^{t(z-1)}$ which is analytic in the entire plane whereas for $p>0$ the generating
function $g_p(z)$ has an essential singularity at $\frac1p$. The expansion at $z_0=1$ yields
$g_p(z)=\sum_{k=0}^\infty (z-1)^kQ_k$ where
$$
Q_k(t,p)=\frac1{(1-p)^k}\sum_{j=1}^kp^{k-j}\frac{t^j}{j!}
\left(\begin{array}{c}k-1\\j-1\end{array}\right)
$$
($Q_0=1$) are the factorial moments. Note that in particular $P_0(0,p)=1$ and $P_r(0,p)=0$ for all
$r\ge1$.

\subsection{Return times patterns}\label{returntimespatterns}

Let $M$ and $m<M$ be given integers (typically $m<\!\!<M$) and let $\tau\in\mathbb{N}$
 be some (large) number. For $r=1,2,3,\dots$ we define the following:\\
{\bf (I)} $G_r(\tau)$: We denote by $G_r(\tau)$ the $r$-vectors
$\vec{v}=(v_1,\dots,v_r)\in\mathbb{Z}^r$ for which $1\leq v_1<v_2<\cdots<v_r\leq \tau$.\\
{\bf (II)} $G_{r,j}(\tau)$: We divide the set $G_r$ into disjoint subsets $G_{r,j}$ where $G_{r,j}$
consists of all $\vec{v}\in G_r$ for which we can find $j$ indices
$i_1,i_2,\dots,i_j\in\{1,2,\dots,r\}$, $i_1=1$, so that $v_k-v_{k-1}\le M$ if
$k\not=i_2,\dots,i_j$ and so that $v_k-v_{k-1}>M$ for all $k=i_2,\dots,i_j$. 

For $\vec{v}\in G_{r,j}$ the values of 
$v_i$ will be identified with returns; returns that occur within less than time $M$ are called
{\em immediate returns} and if the return time is $\ge M$ then we call it a {\em long return}
(i.e.\ if $v_{i+1}-v_i<M$ then we say $v_{i+1}$ is an immediate return and if $v_{i+1}-v_i\ge M$
the we call $v_i$ a long return). 
That means that  $G_{r,j}$ consists of all return time patterns $\vec{v}$ which have $r-j$
immediate returns that are clustered into $j$ blocks of immediate returns and $j-1$ long returns between those blocks. The entries $v_{i_k}$, $k=1,\dots,j$, are the beginnings (heads) of the blocks
(of immediate returns). We assume from now on that all short returns are multiples of $m$. 
(This reflects the periodic structure around periodic points, cf.\ condition~(II) of 
Proposition~\ref{sevastyanov}.)\\
{\bf (III)} $G_{r,j,w}(\tau)$:  For $\vec{v}\in G_{r,j}$ the length of each block is 
$v_{i_{k+1}-1}-v_{i_k}$, $k=1,\dots,j-1$. Consequently 
let us put $w_k=\frac1m(v_k-v_{k-1})$ for the {\em individual overlaps}, for $k\not=i_1,i_2,\dots,i_j$. Then
$\sum_{\ell=i_k+1}^{i_{k+1}-1}w_\ell=\frac1m(v_{i_{k+1}-1}-v_{i_k})$ is the {\em overlap} of the
 $k$th block and
$w=w(\vec{v})=\sum_{k\not=i_1,i_2,\dots,i_j}w_k$ the {\em total overlap} of $\vec{v}$.
We now put
$G_{r,j,w}=\{\vec{v}\in G_{r,j}: w(\vec{v})=w\}$.  ($G_{r,j}=\bigcup_{w}G_{r,j,w}$ is a disjoint
union.)\\
{\bf (IV)} $\Delta(\vec{v})$: For $\vec{v}$ in $G_{r,j}$ we put
$$
\Delta(\vec{v})=\min\left\{v_{i_k}-v_{i_k-1}:\;k=2,\dots,j\right\}
$$
for the minimal distance between the `tail' and the `head' of successive blocks of immediate returns
(or the length of the shortest one of the long gaps).

\subsection{Compound Poisson approximations}

The purpose of this  section is to prove the following result on the approximation of
the compound Poisson distribution.

\begin{proposition}\label{sevastyanov} Let $M, m,\tau$ be as above.
Let $\eta_j$, $j=1,\dots,\tau$, be $0,1$-valued random variables on some $\Omega$
for $\vec{v}\in G_r$ put $\eta_{\vec{v}}=\prod_i\eta_{v_i}$.
Choose $\delta>0$ and define the `rare set' $R_r=\bigcup_{j=1}^rR_{r,j}$, where $R_{r,j}=\{\vec{v}\in
G_{r,j}: \Delta(\vec{v})<\delta\}$.
Let $\mu$ be a probability measure on $\Omega$ which satisfies the following conditions
($C_0$ is a constant):\\
(I) $\mathbb{E}(\eta_j)=\beta$ for all $j=1,\dots,\tau$ (invariance of the measure).\\
(II) Suppose that there are numbers $0<p_-\le p\le p_+$, $\phi\ge0$ so that for all 
$\vec{v}\in G_{r,j,w}\setminus R_{r,j}$
$$
\left|\mathbb{E}(\eta_{\vec{v}})-p^w\beta^j\right|\le C_0\beta^j (p_+^w-p_-^w)+p^w\left((1+\phi)^j-1\right)
$$
if all of the individual overlaps $w_\ell$ are multiples of $m$, and 
$$
\mathbb{E}(\eta_{\vec{v}})=0
$$
if some of the individual overlaps $w_\ell$ are not multiples of $m$.\\
(III) There are some constants $\gamma\ge1, \gamma_1, \gamma_2$ small
(e.g.\ $\gamma(\gamma_1+\gamma_2)<\frac1{12}$), so that for all $r$
$$
\sum_{\vec{v}\in R_r}\mathbb{E}(\eta_{\vec{v}}) \le C_0r\gamma^{r}\sum_{j=2}^r\sum_{s=1}^{j-1}
\left(\begin{array}{c}j-1\\s-1\end{array}\right) \gamma_1^{j-s}\frac{(\tau\beta)^s}{s!}
\left(\begin{array}{c}r-1\\j-1\end{array}\right)\gamma_2^{r-j}.
$$
Let us put $\zeta=\sum_{j=1}^\tau\eta_j$ and $t=(1-p)\tau\beta$.

Then there exists a constant $C_1$ so that for every $t>0$ one has
$$
\left|\mathbb{P}(\zeta=r)-e^{-t}P_r(t,p)\right|
\le C_1(\gamma_1+\delta\beta)t^{r-1}\frac{e^{2r}}{r!}
+C_1\left(p_+^{\frac{M}m}+p_+-p_-+\phi\right)\left\{\begin{array}{lll}
\frac{t^r}{r!}e^{2r+\frac52t}&\mbox{if}&t>\frac12pr\\
(2p)^re^{t\frac{1+2p}{1-4p}}&\mbox{if}&t\le\frac12pr\end{array} \right..
$$
\end{proposition}

\noindent
Note that the constants $\gamma$ and $\gamma_2$ don't enter the final 
estimate in an explicit 
way. The significant quantity here is $\gamma_1$ which typically is $<\!\!<1$
where $\gamma,\gamma_2$ only have to be small enough of order ${\cal O}(1)$.

The choice of $\delta$ is central to the application of this proposition. In the application
however $\phi$ depends on $\delta$ and in fact $\phi(\delta)\rightarrow0+$ as
$\delta\rightarrow\infty$.
Obviously a larger value for $\delta$ increases the error term as one sees from the
 expression, but also a smaller value increases the error term since the 
rare set $R_r$ becomes larger and the `mixing property' in (II) will require larger $p_+$ and
smaller $p_-$, thus again increasing the error estimate. The trick is to optimise $\delta$.

\vspace{3mm}

\noindent {\bf Proof.} We compare the generating function $\psi(z)$ for the process $\zeta$
with the generating function $g_p(z)$ for the compound Poissonian. In part (A) we 
compare their Taylor coefficients at $z=1$ and in part (B) we use Cauchy estimates to 
compare their Taylor coefficients at $z=0$ which then gives us the final result. 

\vspace{3mm}

\noindent {\bf (A)} The coefficients at $z=1$ (factorial moments) of the generating function 
$\psi(z)=\sum_{r=0}^\infty z^r\mathbb{P}(\zeta=r)=\sum_{r=0}^\infty (z-1)^rU_r$
are
$$
 U_r=\sum_{\vec{v}\in G_r}\mathbb{E}(\eta_{\vec{v}}),
 $$
while the coefficients of the generating function  $g_p(z)=\sum_{r=0}^\infty (z-1)^rQ_r(t,p)$ for the
 compound Poisson distribution are
 $$
Q_r(t,p)=\frac1{(1-p)^r}\sum_{j=1}^r\frac{t^j}{j!}
\left(\begin{array}{c}r-1\\j-1\end{array}\right)p^{r-j}
=\sum_{u=r-j}^\infty\sum_{j=1}^rp^u\beta^j\frac{\tau^j}{j!}
\left(\begin{array}{c}u-1\\r-j-1\end{array}\right) \left(\begin{array}{c}r-1\\j-1\end{array}\right),
$$
where $t=(1-p)\beta\tau$. We will now compare the coefficients $Q_r$ to the 
coefficients $U_r$.
 There are three parts to the comparison:
(i) Assumption~(II) is used to compare the terms for which $\vec{v}\in G_r\setminus R_r$;
(ii) Assumptions~(I) and~(III) are used to estimate the total contributions made
by $\vec{v}\in R_r$;
(iii) We have to estimate the contribution to $Q_r$  that correspond to overlaps $u$ 
which do not occur for vectors $\vec{v}$ in $G_r$ and therefore cannot be matched
with terms in the sum that defines $U_r$.

More precisely, we estimate as follows:
$$
\left|U_r-Q_r((1-p)\beta\tau,p)\right| \le
\sum_j\sum_u\sum_{\vec{v}\in G_{r,j,u}\setminus R_j} \left|\mathbb{E}(\eta_{\vec{v}})-p^w\beta^j\right| 
+\sum_{\vec{v}\in R_r}\left(\mathbb{E}(\eta_{\vec{v}})+p^{u(\vec{v})}\beta^j\right) +V(r).
$$
Before we proceed to bound the three terms on the right hand side let us estimate the cardinality of the 
sets $G_{r,j,u}$. (Note that $u\ge r-j$ if $G_{r,j,u}$  is nonempty.)
 Since $G_{r,j,u}$ consists of all $\vec{v}\in G_{r,j}$  that have a total overlap
$u$ (in $j$ blocks of immediate returns)  we get
$$
|G_{r,j,u}|\le\frac{\tau^j}{j!} \left(\begin{array}{c}u-(r-j)+r-j-1\\r-j-1\end{array}\right)
\left(\begin{array}{c}r-1\\j-1\end{array}\right) =\frac{\tau^j}{j!}
\left(\begin{array}{c}u-1\\r-j-1\end{array}\right) \left(\begin{array}{c}r-1\\j-1\end{array}\right)
$$
($j$ blocks positioned `anywhere' on an interval of length $\tau$, $u$ overlaps distributed on
$r-j$ immediate returns and $j$ blocks beginning on any of the $r$ return times).

Now we estimate the three terms in the coefficient comparison as follows:\\
{\bf (i)} The first error term (difference between the dominating terms) is bounded using assumption (II):
\begin{eqnarray*}
\sum_{j=1}^r\sum_{u=r-j}^\infty\sum_{\vec{v}\in G_{r,j,u}\setminus R_r}
\left|\mathbb{E}(\eta_{\vec{v}})-p^w\beta^j\right|
\hspace{-5cm}\\
&\le&\sum_{j=1}^r\sum_{u=r-j}^\infty|G_{r,j,u}|\beta^j\left(p_+^u-p_-^u(1-\phi)\right)\\
&\le&\sum_{j=1}^r\sum_{u=r-j}^\infty\frac{\tau^j}{j!}
\left(\begin{array}{c}u-1\\r-j-1\end{array}\right) \left(\begin{array}{c}r-1\\j-1\end{array}\right)
\beta^j\left(p_+^u-p_-^u(1-\phi)\right)\\
&\le&\sum_{j=1}^r\frac{\tau^j}{j!} \beta^j \left(\begin{array}{c}r-1\\j-1\end{array}\right)
\left(\left(\frac{p_+}{1-p_+}\right)^{r-j}
-\left(\frac{p_-}{1-p_-}\right)^{r-j}(1-\phi)\right)\\
&\le&\frac{c_1q}{(1-p_+)^2} \sum_{j=1}^{r-1}(r-1)\frac{\tau^j\beta^j}{j!}
\left(\begin{array}{c}r-2\\j-1\end{array}\right)
\left(\frac{p_+}{1-p_+}\right)^{r-j-1}+\phi\sum_{j=1}^r\frac{\tau^j}{j!} \beta^j
\left(\begin{array}{c}r-1\\j-1\end{array}\right)
\left(\frac{p_+}{1-p_+}\right)^{r-j}\\
&\le& c_2 q(r-1)Q_{r-1}(t,p_+)+\phi Q_r(t,p_+)
\end{eqnarray*}
(because $\scriptsize(r-j)\left(\begin{array}{c}r-1\\j-1\end{array}\right)
=(r-1)\left(\begin{array}{c}r-2\\j-1\end{array}\right)$), where $t=(1-p)\tau\beta$ and $q=p_+-p_-$.\\
{\bf (ii)} For the second term let us note that $R_r=\bigcup_j R_{r,j}$ where
 $R_{r,j}=\left\{\vec{v}\in G_{r,j}:\;\Delta(\vec{v})<\delta\right\}$. 
 Put $R_r^s$ for those $\vec{v}\in R_r$ where $v_{i+1}-v_i\geq\delta$ for
exactly $s-1$ indices $i_1, i_2, \dots, i_{s-1}$ and put $i_s=v_r$ (obviously $1\le s\le j-1\leq r-1$ and
$i_{s-1}\leq r-1$). 

 To estimate the cardinality of $R_{r,j,u}^s= R_r^s\cap G_{r,j,u}$ let us note that the number of
possibilities of $v_{i_1}<v_{i_2}\cdots<v_{i_s}$ (entrance times for long returns bigger than
$\delta$) is bounded above by $\frac1{s!}\tau^s$ (this is the upper bound for the number of
possibilities to obtain $s-1$ intervals contained in the interval $[1,\tau]$), and each of the
remaining $j-s$ return times less than $\delta$ assume no more than $\delta$ different values.
Since the indices $i_s,\dots,i_{k_s}$ can be picked in
$\scriptsize\left(\begin{array}{c}j-1\\s-1\end{array}\right)$ many ways out of $j$ blocks, we
obtain:
$$
|R^s_{r,j,u}|\leq\left(\begin{array}{c}j-1\\s-1\end{array}\right) \frac{\delta^{j-s}}{s!}\tau^s
\left(\begin{array}{c}r-1\\j-1\end{array}\right)
\left(\begin{array}{c}u-1\\r-j-1\end{array}\right).
$$
To estimate the contribution made by the portion of the sum in the definition of $Q_r$ 
which corresponds to the vectors $\vec{v}\in R_r$ we obtain by summing over $s$:
\begin{eqnarray*}
\sum_{\vec{v}\in R_r}\beta^jp^{w(\vec{v})} &\le&\sum_j\sum_{s=1}^{j-1}\sum_{u=r-j}^\infty
\beta^jp^u|R_{r,j,u}^s|\\
&\leq&\sum_{j=2}^r\sum_{s=1}^{j-1}\left(\begin{array}{c}j-1\\s-1\end{array}\right)
\frac{(\tau\beta)^s}{s!}(\delta\beta)^{j-s} \left(\begin{array}{c}r-1\\j-1\end{array}\right)
\sum_{u=r-j}^\infty\left(\begin{array}{c}u-1\\r-j-1\end{array}\right)
p^u\\
&\leq&\sum_{j=2}^r\sum_{s=1}^{j-1}\left(\begin{array}{c}j-1\\s-1\end{array}\right)
\frac{(\tau\beta)^s}{s!}(\delta\beta)^{j-s} \left(\begin{array}{c}r-1\\j-1\end{array}\right)
\left(\frac{p}{1-p}\right)^{r-j}.
\end{eqnarray*}
The corresponding term for the actual expected values of $\eta_{\vec{v}}$ where $\vec{v}$ 
is in the rare set is bounded by assumption~(II). Hence we obtain
$$
\sum_j\sum_{\vec{v}\in R_{r,j}}(\mathbb{E}(\eta_{\vec{v}})+\beta^jp^{w(\vec{v})} )
 \le S_r,
$$
where
$$
S_r= c_3r\gamma^{r}\sum_{j=2}^r\sum_{s=1}^{j-1}
\left(\begin{array}{c}j-1\\s-1\end{array}\right) \hat\gamma_1^{j-s}\frac{(\tau\beta)^s}{s!}
\left(\begin{array}{c}r-1\\j-1\end{array}\right)\hat\gamma_2^{r-j},
$$
for some $c_3$,  and  $\hat\gamma_1=\gamma_1+\delta\beta$,
$\hat\gamma_2=\gamma_2+\frac{p}{1-p}$.\\
{\bf (iii)} Since the sum for $Q_r$ contains many terms that cannot be paired
with terms in the sum of $U_r$ let us look at those combinations of $r$, $j$, and 
$u$ that do not correspond to vectors in $G_r$. Let us denote by $V_{r,j,u}$ the number
of elements for the values of $r,j,u$ that occur in the representation of $Q_r$ and 
are not in $G_r$. Those overcounts occur for overlaps which
have lengths $\ge M$. Denote by $\tilde{u}_1, \tilde{u}_2,\dots,\tilde{u}_{r-j}$ the individual
overlaps of these $r-j$ fictitious `immediate returns'. Then, if the first intersection is of length
$\tilde{u}_1$ ($\ge\frac{M}m$), then
$$
V_{r,j,u}\le r\sum_{\tilde{u}_1=\frac{M}m}^{u-(r-j)} \left|G_{r,j,u-\tilde{u}_1}\right|,
$$
Hence (for $u>\frac{M}m+(r-j-1)$ and $j<r$)
$$
V_{r,j,u}\le r\frac{\tau^j}{j!}\sum_{\tilde{u}_1=\frac{M}m}^{u-(r-j-1)}
\left(\begin{array}{c}u-\tilde{u}_1-1\\r-j-1\end{array}\right)
\left(\begin{array}{c}r-1\\j-1\end{array}\right)
=r\frac{\tau^j}{j!}\left(\begin{array}{c}r-1\\j-1\end{array}\right)
\left(\begin{array}{c}u-\frac{M}m\\r-j-1\end{array}\right)
$$
where we used the identity
$\sum_{y=b}^a\scriptsize\left(\begin{array}{c}y-1\\b-1\end{array}\right)
=\left(\begin{array}{c}a\\b\end{array}\right)$. Let us put $H_r$ for the set of ficticious vectors
$\vec{v}$ that have individual overlaps $u_i$ ($i\not=i_1,\dots i_j$) that are not
allowed, i.e.\ where at least one overlap is longer than $\frac{M}m$. 
Then we can estimate the contribution one gets
by counting over $H_r$ as follows:
\begin{eqnarray*}
V(r)&=&\sum_{\vec{v}\in H_r}\beta^jp_+^{u(\vec{v})}\\
&=&r\sum_{j=1}^{r-1}\sum_{u=\frac{M}m+r-j}^\infty
V_{r,j,u}\beta^jp_+^u\\
&\le&r\sum_{j=1}^{r-1}\beta^j\frac{\tau^j}{j!} \left(\begin{array}{c}r-1\\j-1\end{array}\right)
\sum_{u=\frac{M}m+r-j-1}^\infty
p_+^u\left(\begin{array}{c}u-\frac{M}m\\r-j-1\end{array}\right)\\
&=&r\sum_{j=1}^{r-1}\beta^j\frac{\tau^j}{j!}
\left(\begin{array}{c}r-1\\j-1\end{array}\right)p_+^{\frac{M}m-1}
\left(\frac{p_+}{1-p_+}\right)^{r-j}
\end{eqnarray*}
where we used the identity $\scriptsize\sum_{u=a}^\infty
p^u\left(\begin{array}{c}u-1\\a-1\end{array}\right) =\left(\frac{p}{1-p}\right)^a$. Hence
$$
V(r)\le rp_+^{\frac{M}m-1}\sum_{j=1}^{r-1}\beta^j\frac{\tau^j}{j!}
\left(\begin{array}{c}r-1\\j-1\end{array}\right) \left(\frac{p_+}{1-p_+}\right)^{r-j}
=rp_+^{\frac{M}m-1} Q_r\left(t,p_+\right).
$$

Combining the three estimates (i), (ii) and (iii) yields
$$
\left|U_r-Q_r(t,p)\right| \le \left(rp^{\frac{M}m-1}+\phi\right)
Q_r(t,p_+)+c_5 q(r-1)Q_{r-1}(t,p_+)+S_r
$$

\vspace{3mm}

\noindent {\bf (B)} The difference $\varphi(z)=\psi(z)-g_p(z)$ between the two generating functions
splits into two parts, $\varphi=\varphi_1+\varphi_2$, which we analyse separately
($\varphi_1$ reflects the estimates of parts~(A-i) and~(A-iii), and $\varphi_2$
reflects the estimate of~(A-ii)):\\
{\bf (i)}  The function $\varphi_1(z)$ is majorised by the power series
$$
c_6\left(p^{\frac{M}m}+q +\phi\right) \sum_r|z-1|^rrQ_r(t,p_+).
$$
The sum over $r$ is equal to $\frac{d}{dw}e^{t\frac{w}{1-p_+-p_+w}}$ (where $w=|z-1|$) which can
be bounded by $4e^{t\frac{w}{1-p_+-p_+w}}$ if $|w|\le\frac34\frac{1-p_+}{p_+}$. 
A Cauchy estimate with $|z|=R$ now yields
$$
E_1=\frac1{r!}\left|\varphi_1^{(r)}(0)\right| \le \frac4{R^r}e^{t\frac{R+1}{1-2p_+-p_+R}}.
$$
If $t$ is large so that $\frac{r}t<\frac1{2p}$ then we can take $R=\frac{r}t$ and obtain for instance
that (assuming $p_+<\frac1{20}$)
$$
E_1 \le c_7\epsilon
\left(\frac{t}r\right)^re^{\frac52(r+t)} \le c_8\epsilon\frac{t^r}{r!}e^{2r+\frac52t}
$$
using Stirling's formula, where $\epsilon=c_6\left(p_+^{\frac{M}m}+
q+\phi\right)$. If $t$ is small so that for instance $\frac{r}t\ge\frac1{2p}$ then we take
$R=\frac1{2p}$ and thus obtain
$E_1=\frac1{r!}\left|\varphi_1^{(r)}(0)\right| \le c_8\epsilon (2p)^re^{t\frac{1+2p}{1-4p}}$.\\
{\bf (ii)} The second error function $\varphi_2(z)$ is majorised by the power series ($t'=\beta\tau$)
\begin{eqnarray*}
\sum_r|z-1|^rS_r&=&\sum_{r=2}^\infty|z-1|^r\gamma^r\sum_{j=2}^r\sum_{s=1}^{j-1}
\left(\begin{array}{c}j-1\\s-1\end{array}\right)\frac{t'^s}{s!}
\left(\begin{array}{c}r-1\\j-1\end{array}\right)
\hat\gamma_1^{j-s}\hat\gamma_2^{r-j}\\
&=&\exp\frac{t'\gamma|z-1|}{1-(\hat\gamma_1+\hat\gamma_2)\gamma|z-1|}-
\exp\frac{t'\gamma|z-1|}{1-\hat\gamma_2\gamma|z-1|}\\
&\le&6\gamma\hat\gamma_1|z-1| \exp\frac{t'\gamma|z-1|}
{1-(\hat\gamma_1+\hat\gamma_2)\gamma|z-1|}.
\end{eqnarray*}
if $|z-1|\gamma(\hat\gamma_1+\hat\gamma_2)$ is small enough (e.g. $\le \frac13$), 
where we have used the identity
$$
\sum_{r=2}^\infty\sum_{j=2}^r\sum_{s=1}^{j-1}
\left(\begin{array}{c}j-1\\s-1\end{array}\right)\frac{x^s}{s!}
\left(\begin{array}{c}r-1\\j-1\end{array}\right)y^{j-s}z^{r-j}
=\exp\frac{x}{1-y-z}-\exp\frac{x}{1-z}
$$
(develop into a Taylor series with variable $x$ and use the identity
$\sum_{k=\ell}^\infty\scriptsize{\left(\begin{array}{c}k-1\\
\ell-1\end{array}\right)}b^{k-\ell}=(1-b)^{-\ell}$). Hence if we put $|z|=R$ ($R>1$) then
$$
E_2=\frac1{r!}\left|\varphi_2^{(r)}(0)\right|\le c_9\gamma\hat\gamma_1\frac{e^{2\gamma Rt'}}{R^{r-1}}
$$
if we assume that $(R+1)\gamma\hat\gamma_1$ is small enough (e.g.\ $<\frac13$). 
If $R=\frac{r}{2t'}$ then we get (using Stirling's formula)
$$
\frac1{r!}\left|\varphi_2^{(r)}(0)\right|\le c_{10}\hat\gamma_12^rr^2t^{r-1}\frac{e^{r}}{r!}
$$
(if $t'$ is close enough to $t$).

Since
$\mathbb{P}(\zeta=r)=\frac1{r!}\psi^{(r)}(0)$ and $e^{-t}P_r(t,p)=g_p^{(r)}(0)$
we get by combining the estimates~(i) and~(ii)
$$
\left|\mathbb{P}(\zeta=r)-e^{-t}P_r(t,p)\right| \le\frac1{r!}\left|\varphi^{(r)}(0)\right|\le E_1+E_2
$$
(and $1+\log2<2$) from which follows the result of the proposition. \qed

\vspace{3mm}

\noindent In the following we will apply this proposition to situations that typically
arise in dynamical systems.  There the stationarity condition~(I) of the proposition
is implied by the invariance of the measure. The random variables $\eta_j$ will be
the indicator function of a cylinder set pulled back under the $j$th iterate of the map. 
Condition~(II) is then implied by the mixing property (see below Definition~\ref{phi.mixing}).
The most difficult condition to satisfy is~(III) because it involves `short range' interaction
over which one has little control and which require more delicate estimates
(see Lemma~\ref{R.small} below).
A simpler version of Proposition~\ref{sevastyanov} is the following corollary ($m=1$)
which is easily deduced by putting $\gamma_2=0, \gamma=1,\gamma_1=\varepsilon$.

\begin{corollary}
 Let $M, \tau$ be as above.
Let $\eta_j$, $j=1,\dots,\tau$, be $0,1$-valued random variables 
and  $\eta_{\vec{v}}=\prod_i\eta_{v_i}$ for $\vec{v}\in G_r$.
For some $\delta>0$ let $R_r=\bigcup_{j=1}^rR_{r,j}$, where $R_{r,j}=\{\vec{v}\in
G_{r,j}: \Delta(\vec{v})<\delta\}$.
Assume $\mu$ be a probability measure on $\Omega$ which satisfies the following conditions
($C_0$ is a constant, $\varepsilon>0$):\\
(I) $\mathbb{E}(\eta_j)=\beta$ for all $j=1,\dots,\tau$ (invariance of the measure);\\
(II) Suppose there is a $p\in(0,1)$ so that for all $\vec{v}\in G_{r,j,w}\setminus R_{r,j}$:
$$
\left|\mathbb{E}(\eta_{\vec{v}})-p^w\beta^j\right|\le \varepsilon p^w;
$$
for all  $r$, $j$ and $w$;\\
(III) 
$$
\sum_{\vec{v}\in R_r}\mathbb{E}(\eta_{\vec{v}}) \le\varepsilon.
$$ 
for all $r=1,2,\dots$.

Then there exists a constant $C_1$ so that for every $t>0$ one has
($\zeta=\sum_{j=1}^\tau\eta_j$ and $t=(1-p)\tau\beta$)
$$
\left|\mathbb{P}(\zeta=r)-e^{-t}P_r(t,p)\right|
\le C_1(\varepsilon+\delta\beta)t^{r-1}\frac{e^{2r}}{r!}
+C_1\left(p^M+\varepsilon\right)\left\{\begin{array}{lll}
\frac{t^r}{r!}e^{2r+\frac52t}&\mbox{if}&t>\frac12pr\\
(2p)^re^{t\frac{1+2p}{1-4p}}&\mbox{if}&t\le\frac12pr\end{array} \right..
$$

\end{corollary}

\section{Measures that are $(\phi,f)$-mixing}
Let $T$ be a map on a space $\Omega$ and $\mu$ a probability measure on $\Omega$. Moreover let
$\cal A$ be a measurable partition of $\Omega$ and denote by ${\cal
A}^n=\bigvee_{j=0}^{n-1}T^{-j}{\cal A}$ its {\em $n$-th join} which also is a measurable partition
of $\Omega$ for every $n\geq1$. The atoms of ${\cal A}^n$ are called {\em $n$-cylinders}. Let us
put ${\cal A}^*=\bigcup_{n=1}^\infty{\cal A}^n$ for the collection of all cylinders in $\Omega$ 
and put $|A|$ for the length of an $n$-cylinder $A\in{\cal A}^*$, i.e.\ $|A|=n$ if $A\in{\cal A}^n$.

We shall assume that $\cal A$ is generating, i.e.\ that the atoms of ${\cal A}^\infty$ are single
points in $\Omega$.

In the following definition we generalise the `retarded strong mixing condition' (see e.g.\
\cite{PS}). We consider mixing dynamical systems in which the function $\phi$ determines the rate
of mixing while the {\em separation function} $f$ specifies a lower bound for the size of the gap
$m$ that is necessary to get the mixing property.

\begin{definition} \label{phi.mixing}
Assume $\mu$ is a $T$-invariant probability measure on $\Omega$ and 
that there are functions $f$ and $\phi$ so that:\\
(i) $f: \mathbb{N}\rightarrow\mathbb{N}_0$
($\mathbb{N}_0=\mathbb{N}\cup\{0\}$) is non-decreasing\\
(ii) $\phi:\mathbb{N}_0\rightarrow\mathbb{R}^+$ is non-increasing.

We say that the dynamical system $(T,\mu)$ is {\em $(\phi,f)$-mixing} if
$$
\left|\mu(U\cap T^{-m-n}V)-\mu(U)\mu(V)\right|\leq\phi(m)\mu(U)\mu(V)
$$
for all $m\geq f(n)$, $n\ge0$, measurable $V$ (in the $\sigma$-algebra generated by ${\cal A}^*$)
and $U$ which are unions of $n$-cylinders.
\end{definition}

\noindent Systems that are $(\phi,f)$-mixing are for instance:
\begin{enumerate}
\item   Classical $\phi$-mixing systems (see, e.g.\
\cite{CGS}): $f=0$. These include equilibrium states for H\"older continuous
potentials on Axiom A systems (which include subshifts of finite type) or on the 
Julia set of hyperbolic rational maps. In this case the partition $\cal A$ is finite.
\item Dispersing billiards \cite{Pen}: $f$ is linear.
\item Equilibrium states for H\"older continuous potentials (that satisfy the 
supremum gap (see section~5) on the Julia set of rational maps 
where  the Julia set contains critical points: $f$ is linear, $\phi$ is exponential.
\item Multidimensional piecewise continuous maps \cite{Pac}: $f$
depends on the individual cylinders ($|{\cal A}|<\infty$).
\end{enumerate}

\vspace{2mm}

\noindent For $r\geq1$ and (large) $\tau\in\mathbb{N}$ let as above $G_r(\tau)$ be the $r$-vectors
$\vec{v}=(v_1,\dots,v_r)\in\mathbb{Z}^r$ for which $1\leq v_1<v_2<\cdots<v_r\leq \tau$. 
 Let $t$ be a positive parameter, $W\subset\Omega$ and put $\tau=[t/\mu(W)]$ be the normalised time. 
 Then the
entries $v_j$ of the vector $\vec{v}\in G_r(\tau)$ are the times at which all the points in
$C_{\vec{v}}=\bigcap_{j=1}^r T^{-v_j}W$ hit the set $W$ during the time interval $[1,\tau]$.
\begin{lemma}\label{product.mixing}
Let $(T,\mu)$ be $(\phi,f)$-mixing.

Then for all $r>1$, $W_i\subset\Omega$ unions
of $n_i$-cylinders, $i=1,\dots,r$ ($n_i\ge1$), and all  `hitting vectors' $\vec{v}\in G_r(\tau)$ 
with return times $v_{i+1}-v_i\geq f(n_i)+n_i$ ($i=1,\dots,r-1$) one has
$$
\left|\frac{\mu\left(\bigcap_{i=1}^rT^{-v_i}W_i\right)} {\prod_{i=1}^r\mu(W_i)}-1\right|
\leq(1+\phi(d(\vec{v},\vec{n})))^r-1,
$$
and $d(\vec{v},\vec{n})=\min_i(v_{i+1}-v_i-n_i)$.
\end{lemma}

\vspace{3mm}

\noindent A consequence of this is that there exists a $0<\eta<1$ so that for all
$\mu(A)\leq\eta^{|A|}$ for all $A\in{\cal A}^*$.

\subsection{Estimate of the rare set}\label{section.return.times}
In this section we provide an estimate for the rare set for general $(\phi,f)$-mixing maps. We
will then use this result in its full strength later to show that the return times distribution at
periodic points is compound Poissonian for rational maps that have critical points. For a `hitting
vector' $\vec{v}\in G_r(\tau)$ ($\tau$ a large integer) we put $C_{\vec{v}}=\bigcap_{k=1}^r
T^{-v_k}W$. Let $\delta\geq f(|W|)$ ($W$ a union of cylinders of the same lengths) then
$$
R_{r,j}(\tau)=\{\vec{v}\in G_{r,j}(\tau): \min_k(v_{i_k+1}-v_{i_k}-|W|)<\delta\},
$$
where the values $v_{i_1},\dots,v_{i_j}$ are the beginnings of the $j$ blocks of
immediate returns (notation as in section~\ref{returntimespatterns}~(II)).

\begin{lemma}\label{R.small} Assume $(T,\mu)$ is $(\phi,f)$-mixing
and assume that there is an $m\in\mathbb{N}$ so that for every $n$ for which $f(n)\le\delta$ there exists an $M<n$ so that $A_n\cap T^{-\ell}A_n\not=\emptyset$ for
$\ell<M$ implies that $\ell$ is a multiple of $m$.

Then there exists a constant $C_2$ so that
for all $n$-cylinders $A_n$:
$$
\sum_{\vec{v}\in R_r} \mu(C_{\vec{v}}) \leq C_2\gamma^{r-1}\sum_{j=2}^r\sum_{s=1}^{j-1}
\left(\begin{array}{c}j-1\\s-1\end{array}\right) (\delta\mu(A_{n'}))^{j-s}\frac{(\tau\mu(A_n))^s}{s!}
\left(\begin{array}{c}r-1\\j-1\end{array}\right)(\gamma\mu(A_{m'}))^{r-j},
$$
where:\\
(i) $n', m'$ ($m'\le n'$) satisfy $f(n')\le M-n'$ and $f(m')\le m-m'$,\\
(ii) $\gamma>1+\phi(\min_i(v_{i+1}-v_i)-n')$,\\
(iii)  $A_{n'}\in{\cal A}^{n'}$, $A_n\subset A_{n'}$,\\
(iv) $A_{m'}\in{\cal A}^{m'}$, $A_{n'}\subset A_{m'}$.
\end{lemma}

\noindent {\bf Proof.} As  in section~(A-ii) of the proof of Proposition~\ref{sevastyanov}, 
put $R_{r,j}^s$ for those $\vec{v}\in R_{r,j}$ for which $v_{i+1}-v_i\geq\delta$ for
$s-1$ indices $i_1, \dots, i_{s-1}$ ($i_s=v_r$, $s\le j-1$).
 We consider two separate cases: (I) $s\geq2$ and (II) $s=1$.

\noindent {\bf (I)} Assume  $s\geq2$ and $i_1, i_2, \dots, i_{s-1}$ be the indices
for which $v_{i_k+1}-v_{i_k}\geq\delta\ge f(n)$ for $k=1,\dots,s-1$. All the other differences are
$\geq M$ and smaller than $\delta$. Let $A_{n'}$ be an $n'$-cylinder so that $A_n\subset A_{n'}$ 
where $n'$ is
so that $f(n')\leq M-n'$. Let $j$ be the number of blocks (i.e.\ $\vec{v}\in G_{r,j}$) and let
$i'_1,\dots,i'_j$ be the beginnings of the `blocks of immediate returns' (clearly $s\le j-1$).
There are $r-j$ immediate short returns of lengths $\in[m,n)$. Let us put
\begin{eqnarray*}
W_{i_k}&=&A_n\;\;\;\; \mbox{\rm for }\;\;k=1,\dots,s,\\
W_{i'_k}&=&A_{n'}\;\;\;\; \mbox{\rm for }\;\;k=1,\dots,j,\\
W_i&=&A_{m'}\cap T^{-m}A_{m'}\cap T^{-2m}A_{m'}\cap \cdots\cap T^{-(u_i-1)m}A_{m'} 
\;\;\;\;\mbox{\rm for all}\;\; i\not\in\{i_k:k\}\cup\{i'_k:k\}
\end{eqnarray*}
where $u_i$ is the overlap for the $i$th return (which is an immediate periodic return). By our
choice of $n'$ we have achieved that $v_{i_k+1}-v_{i_k}\ge\delta\ge f(n)$ and $v_{i+1}-v_i\geq
f(n')$ for $i\not=i_k,\;k=1,\dots,s-1$ and $i\not=i'_k,\;k=1,\dots,j$. By 
Lemma \ref{product.mixing} we obtain
$$
\mu\left(C_{\vec{v}}\right) \leq\mu\left(\bigcap_{i=1}^{r}T^{-v_i}W_i\right)
\leq\alpha_1^{r-1}\prod_{i=1}^{r}\mu(W_i)
 \leq\alpha_1^{r-1}\mu(A_{n'})^{j-s}\mu(A_n)^s(\alpha_2\mu(A_{m'}))^u,
$$
($\alpha_1=1+\min(\phi(\delta-n),\phi(n-n'))$, $\alpha_2=1+\phi(m-m')$)
 where the components of $\vec{n}=(n_1,\dots,n_r)$ are
given by $n_{i_k}=n$ for $k=1,\dots,s$ and $n_i=n'$ for $i\not=i_k,\;k=1,\dots,s$), where 
$u=\sum_i u_i$ is the total overlap.

The cardinality of $R_{r,j,u}^s= R_r^s\cap G_{r,j,u}$ has been estimated in part~(A-ii)
of Proposition~\ref{sevastyanov} to be
$$
|R^s_{r,j,u}|\leq\left(\begin{array}{c}j-1\\s-1\end{array}\right) \frac{\delta^{j-s}}{s!}\tau^s
\left(\begin{array}{c}r-1\\j-1\end{array}\right)
\left(\begin{array}{c}u-1\\r-j-1\end{array}\right),
$$
Therefore
$$
\sum_{\vec{v}\in R^s_{r,j,u}} \mu(C_{\vec{v}})
\leq\alpha_1^{r-1}\left(\begin{array}{c}j-1\\s-1\end{array}\right)
\frac{(\tau\mu(A_n))^s}{s!}(\delta\mu(A_{n'}))^{j-s} \left(\begin{array}{c}r-1\\j-1\end{array}\right)
\left(\begin{array}{c}u-1\\r-j-1\end{array}\right)(\alpha_2\mu(A_{m'}))^u.
$$
\noindent {\bf (II)} If $s=1$ then all returns between blocks are less than $\delta$ for all $k$. In the
same way as above we obtain
$$
\sum_{\vec{v}\in R^1_{r,j,u}} \mu(C_{\vec{v}}) \leq\alpha_1^{r-1}\tau\mu(A_n)(\delta\mu(A_{n'}))^{j-1}
\left(\begin{array}{c}r-1\\j-1\end{array}\right)
\left(\begin{array}{c}u-1\\r-j-1\end{array}\right)(\alpha_2\mu(A_{m'}))^u.
$$

\vspace{2mm}

\noindent Summing over $s$ and using the estimates from (I) and (II) yields
\begin{eqnarray*}
\sum_{\vec{v}\in R_r}\mu(C_{\vec{v}}) &=&\sum_j\sum_{s=1}^{j-1}\sum_{u=r-j}^\infty\sum_{\vec{v}\in
R_{r,j,u}^s}
\mu(C_{\vec{v}})\\
&\leq&\sum_{j=2}^r\alpha_1^{r-1} \sum_{s=1}^{j-1}\left(\begin{array}{c}j-1\\s-1\end{array}\right)
\frac{(\tau\mu(A_n))^s}{s!}(\delta\mu(A_{n'}))^{j-s} \left(\begin{array}{c}r-1\\j-1\end{array}\right)
\sum_{u=r-j}^\infty\left(\begin{array}{c}u-1\\r-j-1\end{array}\right)
(\alpha_2\mu(A_{m'}))^u\\
&\leq&\sum_{j=2}^r\alpha_1^{r-1} \sum_{s=1}^{j-1}\left(\begin{array}{c}j-1\\s-1\end{array}\right)
\frac{(\tau\mu(A_n))^s}{s!}(\delta\mu(A_{n'}))^{j-s} \left(\begin{array}{c}r-1\\j-1\end{array}\right)
\left(\frac{\alpha_2\mu(A_{m'})}{1-\alpha_2\mu(A_{m'})}\right)^{r-j}\\
\end{eqnarray*}
The lemma now follows since 
$\frac{\alpha_2\mu(A_{m'})}{1-\alpha_2\mu(A_{m'})}\le\alpha'\mu(A_{m'})$ with a
$\alpha'$ which is slightly larger than $\alpha_2$. Now we write $\alpha_2$ instead of $\alpha'$. \qed

\vspace{3mm}

\noindent In the case of classical $\phi$-mixing maps (see subsection 3.2 below), when $f$ is zero,
we get the following simpler result. (We simply put $n'=n$ and $m'=m$ which then results in
$A_{n'}=A_n$ and $A_{m'}=A_m$.)

\begin{corollary}\label{R.small.phi-mixing} Assume $(T,\mu)$ is $\phi$-mixing
and assume that there is an $m\in\mathbb{N}$ so that for every $n$ there exists an $M<n$ 
so that $A_n\cap T^{-\ell}A_n\not=\emptyset$ for
$\ell<M$ implies that $\ell$ is a multiple of $m$.

Then there exists a constant $C_3$ so that
for all $n$-cylinders $A_n$:
$$
\sum_{\vec{v}\in R_r} \mu(C_{\vec{v}}) \leq C_3\alpha^{r-1}\sum_{j=2}^r\sum_{s=1}^{j-1}
\left(\begin{array}{c}j-1\\s-1\end{array}\right) (\delta\mu(A_n))^{j-s}\frac{(\tau\mu(A_n))^s}{s!}
\left(\begin{array}{c}r-1\\j-1\end{array}\right)(\alpha\mu(A_m))^{r-j},
$$
where $\alpha=1+\phi(0)$ and $A_m\in{\cal A}^m$ contains $A_n$.
\end{corollary}

\subsection{$\phi$-mixing measures}

We say that the dynamical system $(T,\mu)$ is {\em $\phi$-mixing} if $f$ is identically zero, i.e.
$$
\left|\mu(U\cap T^{-m-n}V)-\mu(U)\mu(V)\right|\leq\phi(m)\mu(U)\mu(V)
$$
for all $m$, measurable $V$ (in the $\sigma$-algebra generated by ${\cal A}^*$) and $U$ which are
unions of cylinders of the same length $n$, for all $n$. The function $\phi$ is assumed to be
monotonically decreasing to zero.

Let $W$ be a set in $\Omega$. Then the entries $v_i$ of the vector $\vec{v}\in G_r(\tau)$ are the
times at which all the points in $C_{\vec{v}}=\bigcap_{i=1}^r T^{-v_i}W$ hit the set $W$ during the
time interval $[1,\tau]$. Following Lemma~\ref{product.mixing} we get that for $n_i$-cylinders
$W_i\subset\Omega$, $i=1,\dots,r$:
\begin{equation}\label{phi.mixing}
\left|\frac{\mu\left(\bigcap_{i=1}^rT^{-v_i}W_i\right)} {\prod_{i=1}^r\mu(W_i)}-1\right|
\leq(1+\phi(d(\vec{v})))^r-1,
\end{equation}
for all `hitting vectors' $\vec{v}\in G_r(\tau)$ with return times $v_{i+1}-v_i\geq n_i$
($i=1,\dots r-1$) where $d(\vec{v})=\min_i(v_{i+1}-v_i-n_i)$.

\subsection{Distribution near periodic points for $\phi$-mixing measures}

\begin{lemma} Let $x$ be a periodic point with minimal period $m$.
If $\mu$ is $\phi$-mixing then the limit
$$
p=\lim_{\ell\rightarrow\infty}\left|\frac1\ell \log\mu(A_{\ell m}(x))\right|
$$
exists.
\end{lemma}

\noindent {\bf Proof.} We show that the quantity inside the logarithms is nearly superadditive.
Let $\Delta$ be an integer so that $\phi(\Delta m)\le\frac12$. Then we have
$$
\left|\log\mu(A_{km+\Delta m+\ell m}(x))\right| \ge\left|\log\mu(A_{km} \cap T^{-km-\Delta
m}A_{\ell m}( T^{-km-\Delta m}x))\right|
$$
and by the mixing property
$$
\mu(A_{km}\cap T^{-km-\Delta m}A_{\ell m}(x)) =\mu(A_{km})\mu(A_{\ell m}(x))(1+{\cal
O}^*(\phi(\Delta m)))
$$
(where ${\cal O}^*$ means that $\left|\frac{{\cal O}^*(\varepsilon)}{\varepsilon}\right|\le1$ for
all $\varepsilon$). If we put $a_j=\left|\log\mu(A_{jm}(x))\right|$, then
$$
a_{k+\Delta+\ell}\ge a_k +a_\ell -\left|\log(1-\phi(\Delta m))\right| \ge a_k+a_\ell-2\phi(\Delta
m)
$$
for all positive integers $k,\ell$. Iterating this inequality yields
\begin{eqnarray*}
\frac{a_{rk+(r-1)\Delta+s}}{(rk+(r-1)\Delta+s)m}&\ge&
\frac{ra_{k}-2(r-1)\phi(\Delta m)}{(rk+(r-1)\Delta+s)m}\\
&\ge&\frac1{1+\frac{\Delta}k+\frac{s}{kr}}\frac{a_{k}}{km}-\frac{2\phi(\Delta m)}{k+\Delta},
\end{eqnarray*}
for positive integers $k,r$ and $s\in[0,k+\Delta)$. If we put $n=kr+(r-1)\Delta+s$, $0\le s\le k+\Delta-1$,
and let $r\rightarrow \infty$ we obtain
$$
\liminf_{n\rightarrow\infty}\frac{a_n}{nm}
\ge\frac1{1+\frac{\Delta}k}\frac{a_{k}}{km}-\frac{2\phi(\Delta m)}{k+\Delta}.
$$
Now let $k\rightarrow\infty$ and we finally get
$$
\liminf_{n\rightarrow\infty}\frac{a_{n}}{nm} \ge\limsup_{k\rightarrow\infty}\frac{a_{k}}{km}
$$
which implies the lemma.
\qed

\vspace{3mm}

\noindent As a consequence of the lemma we see that $p\le \eta^m$ for some
$\eta<1$. (This follows from the fact that $m$-cylinders have measure $\le\eta^m$
for some $\eta<1$ \cite{HV}.) In particular $p$ is always strictly less than $1$.

 In the following we shall assume the stronger property that
$p=\lim_{n\rightarrow\infty}\frac{\mu(A_{n+m}(x))}{\mu(A_n(x))}$. This of course implies the limit
in the lemma, but we are not sure whether the reverse implication is generally true. Also put
$q_n=\sup_{\ell\ge n}\left|\frac{\mu(A_{\ell+m}(x))}{\mu(A_\ell(x))}-p\right|$. For $t>0$ and
integers $n$ we put $\zeta_n^t$ for the counting function 
$\sum_{j=0}^{\tau_n} \chi_{A_n(x))}\circ T^j$ with the observation time
$$
\tau_n=\left[\frac{t}{(1-p)\mu(A_n(x))}\right]
$$
(where $x$ is periodic with minimal period $m$).

\vspace{3mm}

\noindent In order to satisfy the assumptions of Proposition~\ref{sevastyanov} we put
$\gamma=\alpha$,
$\gamma_1=\alpha\delta_n\mu(A_n)$ and $\gamma_2=\alpha\mu(A_m))$
and Corollary~\ref{R.small.phi-mixing}.

\begin{theorem}\label{phi-mixing}
Let $(\mu,\Omega)$ be a $\phi$-mixing measure with partition $\cal A$
(finite or infinite), $x$ a periodic point with minimal period $m$ and $p$ and $q_n$ as above.

Then there exists a constant $C_4$ so that for every $\delta>0$ and every $t>0$ one has
$$
\left|\mathbb{P}(\zeta_n^t=r)-e^{-t}P_r\right|
\le C_4\delta\mu(A_n)t^{r-1}\frac{e^{2r}}{r!}
+C_4\left(p^{\frac{n}m} +q_n+\phi(\delta)\right)
\left\{\begin{array}{lll}
\frac{t^r}{r!}e^{2r+\frac52t}&\mbox{if}&t>\frac12pr\\
(2p)^re^{t\frac{1+2p}{1-4p}}&\mbox{if}&t\le\frac12pr\end{array} \right. ,
$$
\end{theorem}

\noindent {\bf Proof.} We use Proposition~\ref{sevastyanov} and have to verify conditions
(I)--(III). From the definition of $p$ and $q_n$ assumption~(I) is clearly satisfied with $p_\pm =
p\pm q_n$.

To verify condition~(II) let $\vec{v}\in G_{r,j,u}$ and let us look at the measure of
$\mu(C_{\vec{v}})$. By (\ref{phi.mixing}) we have 
($\Delta$ is as defined in section~\ref{returntimespatterns})
$$
\left|\mu(C_{\vec{v}})-\prod_{k=1}^j\mu(D_k)\right|
\le((1+\phi(\Delta(\vec{v})-n))^j-1)\prod_{k=1}^j\mu(D_k),
$$
where $D_k$ is the $k$th block, i.e.\
$$
D_k=\bigcap_{\ell=i_k}^{i_{k+1}-1}T^{-v_\ell}A_n(x).
$$
Since $\mu(D_k)=\mu(A_{n+mu_k}(T^{v_{i_k}}(x)))$ we get by definition of $q_n$
\begin{eqnarray*}
\frac{\mu(D_k)}{\mu(A_n)}
&=&\frac{\mu(A_{n+mu_k})}{\mu(A_n)}\\
&=&\frac{\mu(A_{n+m})}{\mu(A_n)}\frac{\mu(A_{n+2m})}{\mu(A_{n+m})}\cdots
\frac{\mu(A_{n+mu_k})}{\mu(A_{n+m(u_k-1)})}\\
&=&\left(p+{\cal O}(q_n)\right)^{u_k}
\end{eqnarray*}
and therefore
$$
\prod_{k=1}^j\frac{\mu(D_k)}{\mu(A_n)} =\prod_{k=1}^j\left(p+{\cal
O}(q_n)\right)^{u_k}=\left(p+{\cal O}(q_n)\right)^u,
$$
where $u=\sum_{k=1}^ju_k$. Hence
$$
\left|\prod_{k=1}^j\mu(D_k)-p^u\mu(A_n)^j\right| \le\mu(A_n)^j\left((p+q_n)^u-p^u\right)
$$
and consequently
$$
\left|\mu(C_{\vec{v}})-p^u\mu(A_n)^j\right|
\le\mu(A_n)^j\left((p+q_n)^u-p^u+p^u((1+\phi(\Delta(\vec{v})-n))^j-1)\right).
$$
Hence, if $\vec{v}\not\in R_{r,j}$ then we get assumption~(II) with $\gamma=\alpha$,
$p_\pm = p\pm q_n$ (and $\gamma(\gamma_1+\gamma_2)\le\frac1{12}$ if $m,n$ are
not too small). Here we use $M=n-m$.

To verify assumption~(III) we use Corollary~\ref{R.small.phi-mixing} . We obtain
$$
\sum_{\vec{v}\in R_r} \mu(C_{\vec{v}}) \leq C_2\alpha^{r-1}\sum_{j=2}^r\sum_{s=1}^{j-1}
\left(\begin{array}{c}j-1\\s-1\end{array}\right) (\delta\mu(A_n))^{j-s}\frac{(\tau\mu(A_n))^s}{s!}
\left(\begin{array}{c}r-1\\j-1\end{array}\right)(\alpha\mu(A_m))^{r-j},
$$
where $\alpha=1+\phi(0)$. Hence condition~(III) of Proposition~\ref{sevastyanov} is satisfied
with $\gamma_1=\delta\mu(A_n(x))$, $\gamma_2=\alpha\mu(A_m(x))$ and $\beta=\mu(A_n(x))$. \qed

\vspace{2mm}

\noindent Let us note that this result applies to finite as well as infinite partitions $\cal A$.
Since here we focus on the recurrence properties around periodic points we do 
not require the condition $\sum_{A\in{\cal A}}-\mu(A)\log\mu(A)<\infty$ (which is 
necessary in order to get finite entropy or the theorem of Shannon-McMillan-Breiman).

\vspace{4mm}

\noindent {\bf Equilibrium states for Axiom A systems:} Let us now assume that $\mu$ is an
equilibrium state for a H\"older continuous function $f$(with pressure zero) on an Axiom A space
(shift space) which has the finite, generating partition $\cal A$ (see~\cite{Bow}). 
Then $\mu=h\nu$ where $h$ is a normalised eigenfunction for the largest eigenvalue
of the transfer operator and $\nu$ is the associated eigenfunction. In particular $\nu$ is
$e^{-f}$-conformal, i.e.\ if $T$ is one-to-one on a set $A$ then $\nu(TA)=\int_Ae^{-f}\,d\nu(x)$.
If we replace $f$ by $\tilde{f}=f+\log h-\log h\circ T$ then $\mu$ is $e^{-\tilde{f}}$-conformal.
Thus, if $x$ is a periodic point with period $m$, then
$$
\mu(A_n(x))=\mu(T^mA_{n+m}(x))=\int_{A_{n+m}(x)}e^{-\tilde{f}^m(y)}\,d\mu(y)
=\mu(A_{n+m}(x))\tilde{q}_ne^{-\tilde{f}^m(x)},
$$
where $\tilde{q}_n$ is a number that can be estimated by
$$
|\log \tilde{q}_n|\le\var_n\tilde{f}^m \le\const(\var_nf+\var_{n+m}\log h+\var_n\log h)
\le\const\theta^n,
$$
for some $\theta\in(0,1)$ (H\"older exponent). Hence
$$
\frac{\mu(A_{n+m}(x))}{\mu(A_n(x))}=p+q_n,
$$
where $p=e^{f^m(x)}$ and $q_n=p(\tilde{q}_n-1)$ can be estimated by $|q_n|\le p\theta^n\const$. In
particular the limit $\lim_{n\rightarrow\infty}\frac{\mu(A_{n+m}(x))}{\mu(A_n(x))}$ exists and
equals $p$. It is known that $\mu$ is $\phi$-mixing where $\phi(k)=\rho^k$ for some
$\rho\in(\theta,1)$. Let us now apply Proposition 2  and in order to minimise the
term $\epsilon_n=C_1\left(p^{\frac{n}m} + q_n+\rho^{\delta_n}\right)$ we 
choose $\delta_n=\frac{\log\mu(A_n(x))}{\log\rho}$. Then
$\epsilon_n\le\const(p^{\frac{n}m}+n\mu(A_n(x)))$ (again $M=n-m$).

\begin{corollary}
Let $\mu$ be an equilibrium state for a H\"older continuous function on an Axiom A system. Then
there exists a constant $C_5$ so that for all periodic points $x$, $t>0$ and $r=0,1,\dots$ one has
($p$ is as above):
$$
\left|\mathbb{P}(\zeta_n^t=r)-e^{-t}P_r\right|
\le C_5\mu(A_n)|\log\mu(A_n(x))|t^{r-1}\frac{e^{2r}}{r!}
+C_5\left(p^{\frac{n}m}+n\mu(A_n(x))\right) \left\{\begin{array}{lll}
\frac{t^r}{r!}e^{2r+\frac52t}&\mbox{if}&t>\frac12pr\\
(2p)^re^{t\frac{1+2p}{1-4p}}&\mbox{if}&t\le\frac12pr\end{array} \right..
$$
\end{corollary}

\vspace{4mm}

\noindent {\bf Algebraically $\phi$-mixing systems.} If we assume that $\mu$ is $\phi$-mixing
(with respect to the partition $\cal A$) where $\phi(k)={\cal O}(k^{-\kappa})$ for some $\kappa>0$,
then let us note that
$$
p+q_n=\frac{\mu(A_{n+m}(x))}{\mu(A_n(x))} \le
\frac{(1+\phi(0))\mu(A_n(x))\mu(A_m(x))}{\mu(A_n(x))}\le c_1\mu(A_m(x))
$$
implies the very rough estimate $q_n\le\mu(A_m(x))$. With
$\delta_n=\mu(A_m(x))^{-\frac1\kappa}$ one now obtains ($n>>m$)
$$
\epsilon_n \le c_2\left(p^{\frac{n}m}+\mu(A_m(x))+\delta^{-\kappa}\
\right) \le c_3\mu(A_m(x)).
$$

\begin{corollary}
Let $\mu$ is $\phi$-mixing and $\phi(k)\sim k^{-\kappa}$ for some $\kappa>0$. Then there exists a
constant $C_6$ so that for all periodic points $x$, $t>0$ and $r=0,1,\dots$ one has ($p$ is as
above):
$$
\left|\mathbb{P}(\zeta_n^t=r)-e^{-t}P_r\right|
\le C_6\mu(A_n)|\log\mu(A_n(x))|t^{r-1}\frac{e^{2r}}{r!}
+C_6n\mu(A_m(x)) \left\{\begin{array}{lll}
\frac{t^r}{r!}e^{2r+\frac52t}&\mbox{if}&t>\frac12pr\\
(2p)^re^{t\frac{1+2p}{1-4p}}&\mbox{if}&t\le\frac12pr\end{array} \right..
$$
\end{corollary}

\subsection{Example}

\noindent In \cite{KL,Lac} it has been shown that for ergodic systems every possible distribution
can be realised for entry and return times of ergodic systems if the sequence of sets is suitably
chosen. Naturally all settings in which the limiting distributions are shown to be exponential or
Poissonian (in the case of higher returns) have to assume that the target set is a cylinder set
(or a topological ball as in \cite{Pitskel,H2}). Here we show that even if we take cylinder sets
then there are points which do not have a limiting distribution at all.

For simplicity's sake let $\Sigma$ be the full two shift with symbols $0,1$ on which we put the
Bernoulli measure with weights $w,1-w>0$ ($w\not=\frac12$). Let $y=0^\infty$ and $z=1^\infty$ be
the two fixed points under the shift transformation $\sigma$. They have periods $m_1=m_2=1$. The
entry times at $y,z$ are compound Poissonian with the $p$-weights $p_1=w$ and $p_2=1-w$. Put
$\varepsilon=\frac13|p_1-p_2|$ and we will now produce a point $x$ so that the return times
distribution up to some order $r_0$ oscillates between the two compound Poisson distributions.
Choose $n_1$ so that the cylinder $A_{n_1}(y)=A_{n_1}(0^{n_1})$ has the distribution
$$
\left|\mathbb{P}(\zeta_{n_1}^t=r)-e^{-t}P_r(t,p_1)\right|<\frac\varepsilon3
$$
for $t\le t_0$ and $r=1,\dots,r_0$ for some $t_0>0$. Now we choose $n_2>n_1$ so that for the
cylinder $A_{n_2}(0^{n_1}1^{n_2-n_1})$ one has
$$
\left|\mathbb{P}(\zeta_{n_2}^t=r)-e^{-t}P_r(t,p_2)\right|<\frac\varepsilon3
$$
for $t\le t_0$ and $r=1,\dots,r_0$. This can be done because the limiting distribution is
invariant under the shift $\sigma$
(i.e.\ the limiting distribution of the cylinder $A_{n_2}(0^{n_1}1^{n_2-n_1})$ as
$n_2\rightarrow\infty$ is equal to the limiting distribution of the cylinder $A_n(1^\infty)$
as $n\rightarrow\infty$).
 Continuing in this way we find a sequence of integers
$n_1,n_2,n_3,\dots$ so that the distribution of $\zeta_{n_j}^t$ alternates within an error of
$\frac\varepsilon3$ between the distribution $e^{-t}P_r(t,p_1)$ (for odd $j$) and
$e^{-t}P_r(t,p_2)$ (for even $j$) for $t\le t_0$ and $r\le r_0$. Hence the point $x=\bigcap_j
A_{n_j}(0^{n_1}1^{n_2-n_1}\cdots *^{n_j-n_{j-1}-\cdots-n_1})$ ($*$ is $0$ is $j$ is odd and $1$ if
$j$ is even) has no limiting distribution.

Naturally, this construction can be carried out in all $\phi$-mixing systems. Instead of two fixed
points one can also take any finite number of periodic points and then construct a point which
takes turns visiting all of those so that at each visit it stays long enough so that its return
time distribution gets arbitrarily close to the return time distribution of the periodic orbit it
visits.

\section{Return times}\label{section.return.times}

Instead of looking at the probability of a randomly chosen point in the space $\Omega$ to enter a
given set $A$, here we look at the statistics with which points within $A$ return to $A$ again. In
the case of the first entry and return times, these two distributions have for general ergodic
systems been linked in \cite{HLV0}. Higher order entry and return times have been related in
\cite{CK}. It turns out that these distributions are the same only if the first return time is
exponential. Similarly, the number of entry and return times have the same distribution if it is
Poissonian. However, near periodic orbits we get for the return times a distribution which is very
similar, namely it is in the limit given by the following compound Poisson distribution.


Let $p\in(0,1)$. If we define
$$
\hat{P}_r(t,p)=\sum_{j=0}^rp^{r-j}(1-p)^{j+1}\frac{t^j}{j!}
\left(\begin{array}{c}r\\j\end{array}\right)
$$
for $r=1,2,\dots$ and $\hat{P}_0=1-p$ then the generating function for the probabilities
$e^{-t}\hat{P}_r$ is
$$
\hat{g}_p(z)=e^{-t}\sum_{r=0}^\infty z^r\hat{P}_r =\frac{1-p}{1-zp}\,e^{t\frac{z-1}{1-pz}}.
$$
The mean of this distribution is $\frac{t+p}{1-p}$ and the variance is $\frac{t+tp+p}{(1-p)^2}$.
Again note that if $p=0$ then we get the Poisson terms $e^{-t}\hat{P}_r(t,0)=e^{-t}\frac{t^r}{r!}$ and
the generating function $e^{t(z-1)}$ which is analytic in the entire plane whereas for $p>0$ the
generating function $\hat{g}_p(z)$ has an essential singularity at $\frac1p$. The expansion at
$z_0=1$ yields $\hat{g}_p(z)=\sum_{k=0}^\infty (z-1)^k\hat{Q}_k$ where
$$
\hat{Q}_k(t,p)=\frac1{(1-p)^k}\sum_{j=0}^kp^{k-j}\frac{t^j}{j!}
\left(\begin{array}{c}k\\j\end{array}\right)
$$
($\hat{Q}_0=1$) are the factorial moments.

For a set $A$ let us now define the random variable 
$\hat\zeta_A=\chi_A\sum_{j=1}^{\tau_n}\chi_A\circ T^j$ and put $\hat\zeta_n^t=\hat\zeta_{A_n(x)}$ 
where $t=(1-p)\tau_n\mu(A_n(x))$; we also denote
with $\mu_n$ the conditional measure to the cylinder $A_n(x)$.  In a similar way we can now prove
the following result.

\begin{theorem}\label{phi-mixing}
Let $(\mu,\Omega)$ be a $\phi$-mixing measure with partition $\cal A$, $x$ a periodic point with
period $m$ and $p$ and $q_n$ as above.

Then there exists a constant $C_7$ so that for every $\delta>0$ and every $t>0$ one has
$$
\left|\mathbb{P}(\hat\zeta_n^t=r|A_n)-e^{-t}\hat{P}_r\right|
\le C_7n\delta\mu(A_n)t^{r-1}\frac{e^{2r}}{r!}
+C_7n\left(p^{\frac{n}m}+q_n+\phi(\delta)\right)
\left\{\begin{array}{lll}
\frac{t^r}{r!}e^{2r+\frac52t}&\mbox{if}&t>\frac12pr\\
(2p)^re^{t\frac{1+2p}{1-4p}}&\mbox{if}&t\le\frac12pr\end{array} \right.,
$$
where $\tau_n=\frac{t}{(1-p)\mu(A_n(x))}$.
\end{theorem}

\noindent If we compare these error terms to the ones for the entry times, we notice the
additional factor $n$ which comes from satisfying the condition~(I) of
Proposition~\ref{sevastyanov} (cf.\ \cite{HV}).

Let us note that for $r=0$ this result has previously been obtained by Hirata \cite{Hirata1} for
equilibrium states for H\"older continuous function on Axiom A systems. Here however we also get
error estimates:
$$
\left|\mathbb{P}(\hat\zeta_n^t=0|A_n)-(1-p)e^{-t}\right| \le C_6\left(p^{\frac{n}m}+n\mu(A_n(x))\right).
$$
Note that if $p>0$ then $\hat{P}_0(0,p)=1-p$ is strictly less than one and $\hat{P}_r(0,p)=p^r(1-p)$
for $r\ge1$. There is a point mass at $t=0$ which corresponds to immediate returns within the
neighborhood of the periodic point. These are clearly geometrically distributed.

\vspace{3mm}

\noindent {\bf Remark:} By adapting a recent remark of Chamo\^itre and Kupsa \cite{CK}, we proved
in \cite{HLV} under the condition of the existence of the asymptotic distribution of successive
return times that the asymptotic distributions for the entry and return times are related by the
formula ($k=1,\dots$)
$$
D_k(t)=\int_0^t\left(\hat{D}_{k-1}(s)-\hat{D}_k(x)\right)\,ds
$$
where $D_k(t)$ is the limiting distribution $\mathbb{P}(\zeta_n^t=k)$ as $n\rightarrow\infty$, and
$\hat{D}_k(t)=\lim_{n\rightarrow\infty}\mathbb{P}(\hat\zeta_n^t=k)$.

\section{Rational Maps}
Let $T$ be a rational map of degree at least $2$ and $J$ its Julia set. Assume that we executed
appropriate branch cuts on the Riemann sphere so that we can define univalent inverse branches
$S_n$ of $T^n$ on $J$ for all $n\geq1$. Put ${\cal A}^n= \{\varphi(J):\varphi\in S_n\}$
($n$-cylinders). Note that the diameters of the elements in ${\cal A}^n$ go to zero as
$n\rightarrow\infty$. Moreover, ${\cal A}^n$ is not the join of a partition, yet they have all the
properties we require.

Let $f$ be a H\"older continuous function on $J$ so that $P(f)>\sup f$ ($P(f)$ is the pressure of
$f$), let $\mu$ be its unique equilibrium state on $J$ and 
$\zeta_n=\sum_{j=1}^{\tau_n}\chi_{A_n}\circ T^{-j}$ the `counting function' which measures
 the number of times a given point returns to the
$n$-cylinder $A_n$ within the normalised time $\tau_n=[t/\mu(A_n)]$. Although $\mu$ is not a Gibbs
measure we showed in \cite{HV} that for almost every $x$
$$
\mathbb{P}(\zeta_n=r)\rightarrow\frac{t^r}{r!}e^{-t},
$$
as $n\rightarrow\infty$.

\begin{theorem}\label{rational.maps}
Let $T$ be a rational map of degree $\geq2$ and $\mu$ an equilibrium state for H\"older
continuous $f$ (with $P(f)>\sup f$).

Then there exists a $\tilde{\rho}\in(0,1)$ and $C_8$ so that for every periodic point $x\in J$ the
return times are approximately compound Poissonian with the following error terms:
$$
\left|\mathbb{P}(\zeta_n^t=r)-e^{-t}P_r\right|\le
C_8\tilde\rho^nt^{r-1}\frac{e^{3r}}{r!}+C_8\tilde\rho^\frac{n}m \left\{\begin{array}{lll}
\frac{t^r}{r!}\sqrt{r}&\mbox{if}&t>\frac12pr\\
(2tr)^re^{\frac{t}p\frac1{1-2p}}&\mbox{if}&t\le\frac12pr\end{array} \right. ,
$$
where $p=e^{-f^m(x)-mP(f)}$ and $m$ is the minimal period of $x$.
\end{theorem}

\noindent The univalent inverse branches $S_n$ of $T^n$ (with appropriate branch cuts) split into
two categories, namely the uniformly exponentially contracting inverse branches $S_n'$ and the
remaining $S_n''=S_n\setminus S_n'$ for which do not contract uniformly. In \cite{H2} we showed
the following result:
\begin{lemma}\label{product.mixing.rational} (\cite{H2} Lemma 9)
Let $\eta\in(0,1)$. Then there exists a constant $\upsilon>0$ so that for all $r\ge1$ and
$\vec{v}=(v_1,v_2,\dots,v_r)\in G_r$ satisfying $\min_j(v_{j+1}-v_j)\geq(1+\upsilon)n$ (clearly
$r<\frac{\tau_n}{(1+\upsilon)n}$):
$$
\left|\frac{\mu(\bigcap_{j=1}^rT^{-v_j}W_j)}{\prod_{j=1}^r\mu(W_j)}-1\right|\leq \eta^n,
$$
for all sets $W_1,\dots,W_r$ each of which is a union of atoms in ${\cal A}^n$ and for all large
enough $n$.
\end{lemma}

\noindent Let us define the rare set $R_r$: We put $R_r$ for the set all $\vec{v}\in G_r(\tau_n)$ for
which $\min_j(v_{j+1}-v_j)\leq (1+q)n$.
\begin{lemma}\label{rational.periodic}
Let $x\in J$ be a periodic point with (minimal) period $m$. For all large enough $n$ one has that
$A_n(x)\cap T^{-\ell}A_n(x)\not=\emptyset$ for $\ell<n/2$ only if $\ell$ is a multiple of $m$.
\end{lemma}

\noindent {\bf Proof.} Put $n=km+n'$, $0\le n'<n$, and $\phi=\psi^k\cdots\psi^1\phi^{n'}$, where
$\psi^1,\dots,\psi^k\in S_m$, $\phi^{n'}\in S_{n'}$. Since $x\in A$ is periodic with period $m$ we
get that $T^{im}A\cap A\not=\emptyset$ and in particular $x\in T^{im}A$ for all $i=1,\dots,k$.
Since the sets $\psi(J\cap\Omega_m)$ are all disjoint for different $i$, we obtain $\psi^i=\psi^1$
for all $i$. Put $\psi=\psi^1$ and we get $\phi=(\psi)^k\phi^k$ (with $\psi$ concatenated $k$
times).

Now assume that $A\cap T^{-\ell}A\not=\emptyset$ for some $\ell<\frac{n}2$ which is not a multiple
of $m$. Since for some $i$, $im<\ell<(i+1)m$ and $T^{im}A\cap T^{-\ell+im}A\not=\emptyset$, we can
assume that $\ell<m$. Suppose that there are arbitrarily large $n$ so that $\ell<m$ and $V=A\cap
T^{-\ell}A\not=\emptyset$. Similarly as above we put $n=k'\ell +n''$ ($0\le n''<\ell$) and obtain
that $\phi\in S_n$ decomposes as $\phi=(\tilde\psi)^{k'}\tilde\phi^{n''}$ where $\tilde\psi\in
S_\ell$, $\tilde\phi^{n''}\in S_{n''}$.

Now since $(\tilde\psi)^{k'}(J\cap\Omega_\ell)\rightarrow x$ as $k'\rightarrow\infty$, and $x$ is
periodic with period $m$, we see that such $\ell< m$ cannot exist. Hence, for all $n$ large enough
$T^\ell A\cap A\not=\emptyset$ and $\ell<\frac{n}2$ implies that $\ell$ is a multiple of the
period $m$. \qed

\vspace{3mm} \noindent {\bf Proof of Theorem~\ref{rational.maps}.} We are going to verify the
conditions of Proposition~\ref{sevastyanov}. Let
$x\in J$ be periodic with minimal period $m$. Then\\
(I) holds by invariance of the measure $\beta=\mu(A_\varphi)$ for all $j$.\\
(II) Since $\mu=h\nu$ where $h$ is a H\"older continuous density and $\nu$ is $e^{-f}$-conformal
we obtain as before that
$$
\mu(A_n(x))=\mu(T^mA_{n+m}(x))=\int_{A_{n+m}(x)}e^{-\tilde{f}^m(y)}\,d\mu(y)
=\mu(A_{n+m}(x))\tilde{q}_ne^{-\tilde{f}^m(x)},
$$
where we have used that fact that $\mu$ is $e^{-\tilde{f}}$-conformal with respect to the function
$\tilde{f}=f+\log h-\log h\circ T$. The factor $\tilde{q}_n$ satisfies $|\log
\tilde{q}_n|\le\var_n\tilde{f}^m\le\const\theta^n$, for some $\theta\in(0,1)$. Hence
$$
\frac{\mu(A_{n+m}(x))}{\mu(A_n(x))}=p+q_n,
$$
where $p=e^{f^m(x)}$ (independent of $n$) and the error term $q_n=p(\tilde{q}_n-1)$ is bounded as
$|q_n|\le c_1p\theta^n$ for a
constant $c_1$ which is independent of the periodic point $x$.\\
(III) Here we use Lemma~\ref{R.small}. By Lemma~\ref{rational.periodic} we can choose $M=[n/2]$.
Furthermore we set $\delta=(1+\upsilon)n$. According to Lemma~\ref{product.mixing.rational} our
separation function $f$ is given by $f(k)=(1+\upsilon)k$. Hence $n'=[n/(1+\upsilon)]$ and
$m'=[m/(1+\upsilon)]$. Then $A_{n'}$ is the $n'$-cylinder that contains $A_n=A_n(x)$ 
and whose measure is $\mu(A_{n'})\leq\rho^{n/(1+\upsilon)}$. 
Similarly $A_{m'}$ is the $m'$-cylinder that contains $A_m(x)$ and
and whose measure is $\mu(A_{m'})\leq\rho^{m/(1+\upsilon)}$. Let us choose $\tilde\rho<1$ so that
$\tilde\rho>\max\left(\rho^{\frac1{1+\upsilon}},\eta,\vartheta\right)$. Then (for all large enough
$n$)
$$
\sum_{\vec{v}\in R_r} \mu(C_{\vec{v}}) \leq C_2\gamma^{r-1}\sum_{j=2}^r\sum_{s=1}^{j-1}
\left(\begin{array}{c}j-1\\s-1\end{array}\right) \gamma_1^{j-s}\frac{\beta^s}{s!}
\left(\begin{array}{c}r-1\\j-1\end{array}\right)\gamma_2^{r-j},
$$
where $\gamma_1=\delta\mu(A_{n'})\le\tilde\rho^n$, $\gamma_2\le\alpha\mu(A_{m'})\le \tilde\rho^m$,
$\beta=\tau\mu(A_n)$ and by Lemma~\ref{product.mixing.rational} $\alpha=1+\eta^{\delta-n'}$.
Moreover, since $p_+=p+q_n\le \tilde\rho^m$, $\frac\delta\tau\le\tilde\rho^n$, $p_+-p_-\le q_n\le
\rho^n$ and $\phi=\phi(\delta-n')\le\eta^{\delta-n'}\le\tilde\rho^n$ one has
$$
p_+^{\frac{n}m}+p_+-p_-+\phi \le c_1\tilde\rho^\frac{n}m
$$
for some $c_1$. The theorem now follows from Proposition~\ref{sevastyanov}. \qed


\end{document}